\documentclass[10pt,NFSS]{article}

\usepackage[mathscr]{eucal}
\usepackage{amsfonts}
\usepackage{amsmath}
\usepackage{amssymb}
\usepackage{theorem}

\def\Re{\hbox{\rm I\kern-.18em R}}

\newtheorem{Th}{Theorem}

\newtheorem{Prop}{Proposition}

\newtheorem{Cor}{Corollary}

{\theorembodyfont{\rmfamily}}

{\theorembodyfont{\rmfamily}\newtheorem{Rem}{Remark}}

{\theorembodyfont{\rmfamily}\newtheorem{Def}{Definition}}
\mathchardef\d="0064

\allowdisplaybreaks

\begin{document}

\begin{center}
{\LARGE  {\bf{Initial Value Problems of the Sine-Gordon Equation
and Geometric Solutions\\\small Preprint}}}
\end{center}

\begin{center}
{\bf {\Large Magdalena Toda\footnote{Department of Mathematics and
Statistics, Texas Tech University, Lubbock, TX 79409-1042, USA
{\em mtoda@math.ttu.edu}}\\}} \vspace*{.5cm}
\end{center}

{\abstract Recent results using inverse scattering techniques
interpret every solution $\varphi (x,y)$ of the sine-Gordon
equation as a non-linear superposition of solutions along the axes
$x=0$ and $y=0$. This has a well-known geometric interpretation,
namely that every weakly regular surface of Gauss curvature
$K=-1$, in arc length asymptotic line parametrization, is uniquely
determined by the values $\varphi(x,0)$ and $\varphi(0,y)$ of its
coordinate angle along the axes. We introduce a generalized
Weierstrass representation of pseudospherical surfaces that
depends only on these values, and we explicitely construct the
associated family of pseudospherical immersions corresponding to
it.}

\date{}

\

\

{\bf The sine-Gordon equation and initial value problems}

Let $u:D\subset\mathbb{R}^2\longrightarrow \mathbb{R}$ represent a
differentiable function on some open, simply-connected domain $D$.

In [Kri] it had already been shown that every solution $u(x,y)$ of
the sine-Gordon equation
\begin{equation}
u_{xy} = \sin u
\end{equation}
represents ``some type of nonlinear superposition of solutions
$u_1(x,0)$ and $u_2(0,y)$", that is, travelling along different
characteristics. The purpose of this report is {\it to obtain} all
smooth solutions $u(x,y)$ by algebro-geometric methods which
replace the classical ones (such as direct integration, inverse
scattering and numerical integration).

A differentiable solution $\varphi (x,y)$ of (1) represents the
{\it Tchebychev angle} (i.e., angle between arc length asymptotic
coordinate lines) of a weakly regular pseudospherical surface,
measured at the point corresponding to $(x,y)$. By {\it weakly
regular} surface we mean a parametrized surface whose partial
velocity vector fields never vanish, but are allowed to coincide
at a set of points of measure zero. Obviously, at those
singularity points, the parametrization fails to be an immersion.

Thus, every smooth solution $\varphi(x,y)$ of the equation (1)
corresponds to a weakly regular pseudospherical surface. It is
known that every such surface is completely determined by a pair
of arbitrary smooth functions $\alpha(x)$ and $\beta(y)$, such
that $\alpha(x)=\varphi(x,0)$ and $\beta(y)=\varphi(0,y)$. We view
this pair of functions as a {\it pseudospherical analogue of the
Weierstrass representation} from minimal surfaces, and we call it
{\it generalized Weierstrass representation of pseudospherical
surfaces}. We deduced this representation by analogy to a method
presented in [DPW]. Our representation simply turned out to depend
only on the initial values of the Tchebychev angle, $\alpha(x)=
\varphi(x,0)$ and $\beta(y)=\varphi(0,y)$.

The author of this report found this representation in 1998, while
she was a graduate student. At that time, she was not aware of
some outstanding works like [Kri] and [Bo, Ki]. No previous paper
contained a representation for pseudospherical surfaces of type
Weierstrass, and the holomorphic potential of [DPW] that inspired
this approach had only been studied for some harmonic maps (not
for the Lorentz-harmonic maps, like in our case).

However, after it was computed in the spirit of [DPW], this
representation turned out to be characterized by the initial
conditions of a Goursat problem, so we would now like to recall
the following:

\begin{Def}A {\it nonlinear hyperbolic system of equations}
is a system of partial differential equations for functions $ U,V
: D \rightarrow \mathbb{R}$, where $D:=[0,x_0]\times[0,y_0]$:
\begin{equation}
V_x = f(U,V),\qquad U_y = g(U,V),
\end{equation}
with smooth given functions $f, g: \mathbb{R}^2\rightarrow
\mathbb{R}$. We will call {\it initial value problem for a
nonlinear hyperbolic system} the problem consisting of equations
(2), together with the initial conditions

\begin{equation}
U(x,0)= U_0(x), V(0,y)=V_0(y),
\end{equation}

for $(x,y)\in D$. The
functions $U_0:[0,x_0]\rightarrow \mathbb{R}$ and
$V_0:[0,y_0]\rightarrow \mathbb{R}$ are also assumed to be smooth.
\end{Def}

\begin{Prop} (see [Bo3])
The initial value problem for a nonlinear hyperbolic system has a
unique classical solution.
\end{Prop}
For details, see [Bo3], Theorem 1 and its corollary.

Any nonlinear equation of hyperbolic type can be brought to the
form (1), by substitutions of type $U=U(u,u_x), V=V(u,u_y).$

For the particular case of the sine-Gordon equation, one
introduces the independent variables $U=u, V=u_x$ which satisfy a
system of the form (1), namely $U_x =V, V_y = \sin U,$ with
initial conditions (3).

We provide a method of obtaining solutions to such a problem, by
solving a simplified ODE system, followed by a loop group
factorization.

Since many readers are not familiar with this type of
computations, we provided complete arguments for all of our
techniques and results, while also striving for brevity.

\

{\bf Geometric solutions to the sine-Gordon equation}

We begin our study of surfaces with constant negative Gaussian
curvature $K=-1$, called {\it pseudospherical surfaces}, or {\it
K-surfaces}. We recall that all such surfaces are described by a
sine-Gordon equation, with a corresponding Lax system. Let $M$ be
the image of $D=[0,x_0]\times[0,y_0]$ through the differentiable
map $\psi:D\rightarrow \mathbb{R}^3$, where $\psi$ represents {\it
a weakly regular asymptotic line parametrization} (i.e., such that
the coordinate lines are asymptotic lines, and partial velocities
never vanish, so we can assume them to be unitary). An arc length
asymptotic line parametrization is also called {\it Tchebychev
parametrization}.

Let $\varphi$ represent the angle between the asymptotic lines. We
will call it {\it Tchebychev angle}. Singularities of weakly
regular surfaces occur at those values $(x,y)$ where this angle,
$\varphi (x,y)$ equals $0$ or $\pi$. The first fundamental form is
([Ei], [Bo2]):$$ {\rm I}=|\d\psi|^2=\d x^2+2\cos\varphi\d x\d y+\d
y^2.$$ Let $N$ define the normal vector field to the surface (or
Gauss map). Remark that the unit vector field $N$ is orthogonal to
$\psi_x$, $\psi_y$, $\psi_{xx}$, $\psi_{yy}$.

The following obvious result is due to Lie (around 1870) and is of
crucial importance (see also [Bo2]):

\begin{Th} Every pseudospherical surface has a
one-parameter family of deformations preserving {\rm the second
fundamental form}$$ {\rm II}=\sin\varphi \cdot \d x\d y,$$ {\rm
the Gaussian curvature} $K=-1$, and {\rm the angle} $\varphi$
between the asymptotic lines. The deformation is generated by the
transformation $x\mapsto x^*=\lambda^{-1}x$ and $y\mapsto
y*=\lambda y$, $\lambda>0$. (Angle is preserved in the sense that
$\varphi^*(x^*,y^*)=\varphi(x(x^*), y(y^*))$).
\end{Th}

We will refer to this simple change of coordinates as the {\it
Lie-Lorentz transformation}. Lie-Lorentz transformations of a
certain pseudospherical immersion represent its {\it associated
family}, denoted as $\psi^{\lambda }:D\to \Bbb R^{3}$. In order to
define an orthonormal frame on the surface, we consider the
so-called curvature line coordinates, defined by $u_1=x+y$,
$u_2=x-y.$ Partial velocities with respect to $u_1$ and $u_2$ are
orthogonal. This reparametrization diagonalizes both the first and
the second fundamental form. The eigenvectors of the shape
operator are the orthonormal vectors $e_1$ and $e_2$, called
principal directions.
\begin{Def}
For any (weakly regular) pseudospherical immersion
$\psi:D\rightarrow \mathbb{R}^3$, we identify the {\it orthonormal
standard frame} $F=\{\psi, e_1, e_2, N\}$ with the $\rm
SO(3)$-valued function $(e_1, e_2, N)$ defined at every point of
the surface.

We will generically call {\it rotated frame}  $F_{\theta}$ the
frame obtained by rotating the standard frame $F$ by the angle
$\theta(x,y)$ around $N$, in the tangent plane.

In particular for $\theta=\varphi/2$, where $\varphi(x,y)$ is the
Tchebychev angle between the asymptotic directions, the resulting
frame is denoted ${\cal U}:=F_{\varphi/2}$ and is called the {\it
normalized frame} associated with the standard frame $F$ (see
[Wu1], p.18). Expressed in Tchebychev coordinates, the normalized
frame $\cal U$ is oriented just like $F$, and consists of $\psi$,
$\psi_x$, a unit vector orthogonal to $\psi_x$, $\psi_x^\top$, and
the unit normal $N$.

Finally, we will call {\it extended normalized frame} the
normalized frame ${\cal U}^{\lambda}={\cal U} (x,y,\lambda)$
corresponding to the immersion $\psi^{\lambda}$, obtained via
Lie-Lorentz transformation of coordinates from the immersion
$\psi$.
\end{Def}

It is convenient to use $2\times 2$ matrices instead of $3\times
3$ ones. Therefore we recall the Pauli matrices
\begin{equation}
\sigma_1=\begin{pmatrix} 0&1\\1&0\end{pmatrix},\quad
\sigma_2=\begin{pmatrix} 0&-i\\i&0\end{pmatrix},\quad
\sigma_3=\begin{pmatrix} 1&0\\0&-1\end{pmatrix}.
\end{equation}
We identify the $\rm SO(3)$-valued extended normalized frame
${\cal U}^{\lambda}$ with the $\rm SU(2)$-valued {\it function}
$\cal U$ defined on the same domain $D$, with the initial
condition ${\cal U}(0,0,\lambda)=I$, via the spinor
correspondences between $e_k$ ($k=1,2,3$) and matrices ${\cal
U}\cdot i \sigma_k\cdot{\cal U}^{-1}$. We have the following (see
[TU], [Kri], [Bo2], [Bo3]):

\begin{Th}
The extended normalized frame $\cal U^\lambda$ is a $\rm
SU(2)$-valued function of $\lambda>0$, which satisfies the Lax
differential system
\begin{equation}
\partial_x {\cal U^\lambda}={\cal U}^\lambda\cdot{\cal A},\qquad
\partial_y {\cal U^\lambda}={\cal U}^\lambda\cdot{\cal B},
\end{equation}
where
\begin{equation}
{\cal A}=\frac{i}{2}\begin{pmatrix}
\varphi_x&-\lambda\\
-\lambda&-\varphi_x
\end{pmatrix}, \qquad
{\cal B}=\frac{i}{2}{\lambda^{-1}}\begin{pmatrix}
0&e^{-i\varphi}\\
e^{i\varphi}&0
\end{pmatrix}
\end{equation}
The compatibility condition for the system is ${\cal A}_y-{\cal
B}_x-[{\cal A},{\cal B}]=0,$ which can be rewritten as
$\varphi_{xy}=\sin{\varphi}$.

Conversely, given a smooth solution $\varphi(x,y)$ of the
sine-Gordon equation, there exists a unique solution ${\cal
U}(x,y,\lambda)$ of the Lax system. Moreover, this solution is
real analytic in $\lambda$.
\end{Th}

\

{\bf Harmonic maps and the generalized Weierstrass representation}

For a complete characterization of harmonicity in the context of
pseudospherical surfaces, we recommend [Do, St]. Let us remark
that the wave equation $u_{xy}=0$ over the $xy$-plane can be
understood as {\it harmonicity condition} with respect to the
Lorentz metric $dx \cdot dy$. A well-known fact is the following:
if $M$ is a weakly regular surface with $K<0$, then $M$,
considered with its second fundamental form ${\rm II}$ as a
metric, represents a {\it Lorentzian} 2-manifold $(M,{\rm II})$.
The Gauss map $N:(M,{\rm II})\to S^2$ is Lorentz-harmonic (i.e., $
N_{xy} = \rho\cdot N$, where $\rho$ is a certain real-valued
function) {\it iff} the curvature $K<0$ is constant.

It is also well-known that if $M=(D,\psi)$ is, as usual, a
pseudospherical surface given by a Tchebychev immersion
$\psi:D\to\Bbb R^3$, then the frame ${\cal U}:D\to{\rm SU}(2)$
represents a lift of the Gauss map of $N:D\to S^2$, via the
canonical projection relative to the base point $e_3$, namely
$\pi:\rm SU(2)\to S^2\cong \rm SU(2)/S^1$. From this lifting, it
follows (see, for example, [Bo 2]) that the maps $N$ and ${\cal
U}$ are related by the identification $ N\equiv{\cal U}\cdot i
\sigma_3 \cdot {\cal U}^{-1}.$

A very important result obtained by A. Sym ([Sy]) allows us to
obtain the immersion (up to a rigid motion), once we have the
expression of the extended frame. This is presented in several
papers, (e.g. [1, Me, St]):

\begin{Th}
Starting from a given solution $\varphi (x,y)$ of the sine-Gordon
equation, let us consider the initial value problem of the Lax
system with the initial condition ${\cal U}(0,0,\lambda)={\cal
U}_0$. Let $\cal U(\lambda)$ be the solution to this initial value
problem. Then $\cal U(\lambda)$ represents the extended frame
corresponding to the Tchebychev immersion $ \psi^{\lambda} =
{\frac{d}{dt}}{\cal U}^{\lambda}\cdot ({\cal U}^{\lambda})^{-1},$
where $\lambda = e^t$.
\end{Th}

Recall that the Lie algebras ($\rm su(2), [ \cdot, \cdot]$) and
($R^3, \times$) are isomorphic, and so $\psi^{\lambda}$ can be
written immediately in the `classical' way as an immersion in
$\Bbb R^3$.

\

By this result, once we have the extended frame, we can
reconstruct the surface. Since the frame is just a lift $\cal U$
of the Gauss map $N$, we infer that we could reconstruct
everything starting from the Gauss map. However, there is a
freedom in the frame given by a gauge action. Namely, let us act
on the extended normalized frame $\cal U$ via a rotation matrix
$\cal R$. The result is called {\it gauged frame} ${\hat{\cal
U}}$:

\begin{equation}
\hat{\cal U}={\cal R}(0,0)^{-1}\cdot{\cal U}\cdot{\cal R}.
\end{equation}

It will be convenient for our purposes to fix a base point $x_0\in
D,$ e.g. $x_0=(0,0)$, and impose $ {\cal U}(x_0,\lambda)=I$. We
will use this assumption from now on. Also note that the
orthonormal frame $F^{\lambda}$ represents a gauged frame of the
normalized frame $\cal U^{\lambda}$, via a rotation $\cal R$ of
angle $\theta=-\varphi/2$. We have the following consequence of
Theorem 3:

\begin{Cor}
If $F^{\lambda}$ represents the orthonormal frame corresponding to
the associate family of immersions $\psi^{\lambda}$,  then

$\psi^{\lambda}={\cal R}^{-1}({\frac{d}{dt}}{F}^{\lambda}
(F^{\lambda})^{-1}){\cal R},$ where $\lambda = e^t$ and $\cal R$
is the rotation of angle $-\varphi(x,y)/2$.
\end{Cor}
Let us introduce the Cartan connection $\omega^{\lambda}:=-({\cal
U^{\lambda}})^{-1}\d {\cal U^{\lambda}}={\cal A} \ dx + {\cal B} \
dy$, with $\cal A$ and $\cal B$ given by formulas (6). That is,

\begin{equation}
{\omega^{\lambda}}=\frac{i}{2}\begin{pmatrix}
\varphi_x&-\lambda\\
-\lambda&-\varphi_x
\end{pmatrix} \ dx + \frac{i}{2}{\lambda^{-1}}\begin{pmatrix}
0&e^{-i\varphi}\\
e^{i\varphi}&0
\end{pmatrix} \ dy
\end{equation}

Obviously, $\omega^{\lambda}$ represents a ${\Lambda\rm
su}(2)$-valued form, and then it decomposes into a diagonal,
respectively off-diagonal part as $\omega^{\lambda} = \omega_0 +
\omega_1 $, according to the Cartan decomposition of ${\rm
su}(2)$.

The following is a well known result (see [Me, St, 1] and [Me, St,
2]):

\begin{Prop}
There is a one-to-one correspondence between the space of Lorentz
harmonic maps from $D$ to $S^2$ and the equivalence classes of
admissible connections, under the action of the gauge action
introduced above. Moreover, every admissible connection $\omega$
corresponds to its associated loop $\omega^{\lambda}$ satisfying
the flatness condition
\begin{equation}
\d\omega^\lambda+\omega^\lambda\wedge\omega^\lambda=0.
\end{equation}
\end{Prop}

Let further $\omega_0=\omega^{'}_0 + \omega^{''}_0$ and
$\omega_1=\lambda^{-1} \omega^{'}_1 + \lambda \omega^{''}_1$ be
the usual splittings into (1,0) and, respectively, (0,1)-forms,
that is:

\begin{equation}
\omega^{'}_0=\frac{i}{2}\begin{pmatrix}
\varphi_x & 0\\
0 & -\varphi_x
\end{pmatrix} \ dx,\\\qquad
\omega^{''}_0={\textsc{O}},\\\qquad
\omega^{'}_1=\frac{i}{2}\begin{pmatrix}
0 & e^{-i\varphi}\\
e^{i\varphi} & 0
\end{pmatrix} \ dy,\\\qquad
\omega^{''}_1=\frac{i}{2}\begin{pmatrix}
0 & -1\\
-1 & 0
\end{pmatrix} \ dx.
\end{equation}

In this context, we now introduce the twisted loop algebra of
those Laurent polynomials in $\lambda$ with coefficients in $\rm
su(2)$ that are fixed under the $Ad (\sigma_3)$-automorphism, that
is, $$\Lambda{\rm su}(2)^{\rm alg}=\{X: \mathbb{R}_*\rightarrow
\rm su(2);\ X(-\lambda)=\sigma_3 \cdot X(\lambda)\cdot \sigma_3
\}.$$ It will be convenient to use a certain Banach completion of
this algebra. For this purpose, consider the Wiener algebra ${\cal
G}$ that consists of all Laurent series of parameter $\lambda$
with complex-valued coefficients, ${X(\lambda)}=\sum_{k\in\Bbb
Z}X_k\cdot \lambda^k$, with the property that $\sum_{k\in\Bbb
Z}|X_k|<\infty$. We define $\|X(\lambda)\|=\sum_{k\in\Bbb
Z}|X_k|$. Its is well known that this Wiener algebra $\cal G$ is a
Banach algebra relative to this norm, and it consists of
continuous functions. For a matrix $A(\lambda)\in {\rm su}(2,
{\cal G})$, whose entries are elements of $\cal G$, we consider
the norm $\quad\|A\|= {\sum_{i,j=1,2}\|A_{ij}\|},$ where $A_{ij}$
denotes the $(i,j)$-entry of $A$. It can be checked by a direct
computation that $ \quad\|A B\|\leq \|A\| \cdot \|B\|$ and
$\|I\|=1.$ We denote by $${\Lambda{\rm su}}(2):=(\Lambda{\rm
su}(2)^{\rm alg},\|\cdot\|)$$ the completion of $\Lambda{\rm
su}(2)^{\rm alg}$ with respect to this norm. Let us also introduce
the twisted loop group $${\Lambda{\rm SU}(2)}:=\{g\in{\rm SU}(2);\
\sigma_3 g(\lambda)\sigma_3=g(-\lambda)\}.$$ It is well-known that
$\Lambda{\rm SU}(2)$ is a Banach Lie group with Lie algebra ${\rm
Lie}\,\Lambda{\rm SU}(2)=\Lambda{\rm su}(2)$. The twisting
($Ad(\sigma_3)$ invariance) condition on loop algebra $\Lambda{\rm
su}(2)^{\rm alg}$ can be replaced by the following characteristic
property: in spinor representation, the diagonal part is an even
function $\lambda$, while the off-diagonal part is an odd function
of $\lambda$. In order to carry out the construction method of
pseudospherical surfaces, we introduce the following subalgebras
of $\Lambda{\rm su(2)}$:
\begin{equation}
\Lambda^+{\rm su}(2)=\{X(\lambda);\ X(\lambda)\text{
contains only non-negative powers of }\lambda\}\\
\end{equation}
\begin{equation}
\Lambda^-{\rm su}(2)=\{X(\lambda);\ X(\lambda)\text{
contains only non-positive powers of }\lambda\}\\
\end{equation}
\begin{equation}
\Lambda^-_*{\rm su}(2)=\{X(\lambda);\ X(\infty )=0\}
\end{equation}
The connected Banach loop groups whose Lie algebras are described
by definitions above are denoted, respectively, $\Lambda^+{\rm
SU}(2)$, $\Lambda^-{\rm SU}(2)$ and $\Lambda^-_*{\rm SU}(2)$.

In order to obtain the generalized Weierstrass representation of
pseudospherical surfaces, we need to use the following adapted
factorization (introduced in [To2]):

\begin{Th}{\rm (splitting of Birkhoff type, for real parameter
$\lambda$)}

Let $\tilde\Lambda \rm SU(2)$ be the subset of $\Lambda \rm SU(2)$
whose elements, as maps defined on $\Bbb R_+$, admit an analytic
extension to $\Bbb C_*$. It is easy to see that $\tilde\Lambda \rm
SU(2)$ is a subgroup of $\Lambda \rm SU(2)$. Then the
multiplication map $\tilde\Lambda_*^- \rm SU(2)
\times\tilde\Lambda^+ \rm SU(2) \to\tilde\Lambda \rm SU(2)$
represents a diffeomorphism onto the open and dense subset
$\tilde\Lambda_*^- \rm SU(2) \cdot\tilde\Lambda^+ \rm SU(2)$,
called the ``big cell". In particular, if $g\in\tilde\Lambda \rm
SU(2)$ is contained in the big cell, then $g$ has a unique
decomposition$$ g=g_-g_+$$ where $g_-\in\tilde\Lambda_*^- \rm
SU(2)$ and $g_+\in\tilde\Lambda^+ \rm SU(2)$. The analogous result
holds for the multiplication map $\tilde\Lambda_*^+ \rm
SU(2)\times\tilde\Lambda^- \rm SU(2)\to\tilde\Lambda \rm SU(2)$.
\vskip2pt
\end{Th}

This represents a ``linearized" version of the classical {\it
Birkhoff loop group factorization} from [Pr, Se] (where the
splitting was introduced and proved for smooth loops on the unit
circle $S^1$). Note that in [To2], the above theorem was
formulated for $\rm SO(3,\mathbb{R})$, instead of $\rm SU(2)$.
There it was shown that the `Birkhoff' splitting works for
$\lambda$ on {\it any straight-line of the complex plane}.

The first type of Birkhoff factorization, performed away from a
singular set $S_1\subset D$, allows us to split the extended
moving frame ${\cal U}^\lambda:D\to{\rm SU}(2)$ into two parts.
Recall that the first factor of this splitting is of the form
$g_-=I+\lambda^{-1}g_{-1}+\lambda^{-2}g_{-2}+\cdots$, while the
second factor of the splitting is of the form $g_+=g_0+\lambda
g_1+\lambda^2g_2+\cdots$, respectively. Since the ``big cell" is
open and ${\cal U}^\lambda:D\to{\rm SU}(2)$ is continuous, the set
$$ {\tilde D_1}=\{(x,y)\text{ ; }{\cal
U}^{\lambda}(x,y)\text{ belongs to the ``big cell"}\}$$ is open.
Note that $(0,0)\in\tilde D_1$. Let $S_1 = D -\tilde D_1$ denote
the ``singular" set. We have just shown that $S_1$ is closed and
$(0,0)$ is not an element of the set $S_1$. Similarly, we have
$S_2$ and $\tilde D_2$ for the second splitting.

We can perform the two splittings on the extended frame $\cal
U^\lambda$, independently.

Let ${\cal U}={\cal U}^\lambda$ be the extended normalized moving
frame of a pseudospherical surface and let $(x,y)\in
D\setminus(S_1\cup S_2)$. Then, for some uniquely determined
$V_+\in\Lambda^+{\rm SU}(2)$, $V_-\in\Lambda^-{\rm SU}(2)$ and
${\cal U}_-\in\Lambda_*^-{\rm SU}(2)$, ${\cal
U}_+\in\Lambda_*^+{\rm SU}(2)$, ${\cal U}$ can be written as
\begin{equation}
{\cal U}={\cal U}_+\cdot V_-={\cal U}_-\cdot V_+.
\end{equation}
Here ${\cal U}_-$ is an element of the form ${\cal
U}_-=I+\lambda^{-1}{\cal U}_{-1}+\lambda^{-2}{\cal
U}_{-2}+\cdots$, while $V_+$ is an element of the form
$V_+=V_0+\lambda V_1+\lambda^2V_2+\cdots$, respectively. Analogous
expressions can be written for ${\cal U}_+$ and $V_-$,
respectively. We will show that, starting from data of type
Weierstrass, called normalized potentials $\eta ^x$ and $\eta ^y$,
one can obtain the factors ${\cal U}_+$ and ${\cal U}_-$ as
solutions of a simplified ODE system. These two factors represent
the genetic material necessary and sufficient to recreate the
frame and then the immersed surface via the Sym formula.
\begin{Th}
Let ${\cal U}={\cal U}^\lambda$, ${\cal U}_+$ and ${\cal U}_-$ be
as above. Then the following systems of differential equations are
satisfied:
\begin{equation}
({\cal U}_+)^{-1}\cdot{\partial_x{\cal U}_+} =-\lambda\cdot
{\frac{i}{2}}\cdot\begin{pmatrix}
0 & e^{i(\varphi(0,0)-\varphi(x,0))}\\
e^{-i(\varphi(0,0)-\varphi(x,0))} & 0
\end{pmatrix}
\end{equation}
with initial condition ${\cal U}_+(x=0)=I$

and
\begin{equation}
({\cal U}_-)^{-1}\cdot{\partial_y{\cal
U}_-}=\lambda^{-1}\cdot\frac{i}{2} \cdot{\begin{pmatrix}
0 & e^{-i\varphi(0,y)}\\
e^{i\varphi(0,y)} & 0
\end{pmatrix}},
\end{equation}
with initial condition ${\cal U}_-(y=0)=I$.

Moreover, ${\cal U}_+$ does not depend on $y$ and ${\cal U}_-$
does not depend on $x$.
\end{Th}

In some other words, ${\cal U}_+$ and ${\cal U}_-$ are solutions
of some first order systems of differential equations in $x$ and
$y$, respectively.

{\it Proof}. We will prove the first statement. Proving the other
statement is straightforward.

The first Birkhoff splitting implies $ {\cal U}_+={\cal U}\cdot
V_-^{-1}$, which after differentiation gives
\begin{equation}
\d{\cal U}_+=\d{\cal U}\cdot V_-^{-1}-{\cal U}\cdot
V_-^{-1}\cdot\d V_-\cdot V_-^{-1},
\end{equation}
\begin{equation}
{\cal U}_+^{-1}\d{\cal U}_+=V_-({\cal U}^{-1}\d{\cal
U})V_-^{-1}-\d V_-\cdot V_-^{-1}.
\end{equation}
The last equality can also be written as
\begin{equation}
{\cal U}_+^{-1}\d{\cal U}_+=V_-({\cal A} \ dx + {\cal B} \ dy)
V_-^{-1}-\d V_-\cdot V_-^{-1}.
\end{equation}
We will use the Lax equations. In the last equality, we compare
the coefficient of $\d y$ on the left-hand side with the
coefficient of $\d y$ on the right-hand side. The left-hand side
clearly contains only positive powers of $\lambda$, while the
coefficient of $\d y$ on the right-hand side contains non-positive
powers of $\lambda $ only. Thus, ${\cal U}_+$ depends exclusively
on $x$.

Let us now consider the coefficient of $\d x$ in the same
equality. The left-hand side contains only positive powers of
$\lambda$, while the one on the right-hand side, due to the
$\lambda$-dependence of $\cal A$, contains one term in $\lambda$
and no terms in $\lambda^k$, with $k>1$. Next, we can restrict to
a sufficiently small interval around $(0,0)$ on the line $y=0$.
Let now $
V_-=\tilde{V}_0+\lambda^{-1}\tilde{V}_1+\lambda^{-2}\tilde{V}_2+\cdots=\tilde{V}
_0\cdot T_-$, with $T_-\in\Lambda_*^-{\rm SU}(2)$. But since
${\cal U}_+^{-1}(x)\cdot{{\cal U}_+}'(x)$ contains only positive
powers of $\lambda$, we conclude that ${\cal
U}_+^{-1}(x)\cdot{{\cal U}_+}'(x) \d
x=\tilde{V}_0(x,0)\cdot{\omega^{''}_1}\cdot\tilde{V}_0(x,0)^{-1}$,
where $\omega^{''}_1$ is the one from (10). let us now denote
$\tilde{V}_0(x,0) := V_0 $. In order to determine the matrix
$V_0$, one needs to compare the coefficients of the power
$\lambda^0$ in the same equality. As we pointed out, the left-hand
side has positive powers of $\lambda$ only, while the $x$-part of
right-hand side only contains $-V_0\cdot\beta_{0}\cdot V_0^{-1}-\d
V_0\cdot {V_0}^{-1}$ as the only term that does not depend on
$\lambda$, where we denoted
$\beta_0=\omega_0'(x,0)=\frac{i}{2}\begin{pmatrix}
\varphi_x(x,0) & 0\\
0 & -\varphi_x(x,0)
\end{pmatrix} \d x$.
Thus, $V_0$ is a solution to $\d V_0=- V_0 \cdot {\beta_0}$. The
solution $V_0$ of the system must take into account that ${\cal
U}(0,0,\lambda)=I$. Thus $V_0(x)=e^{\theta(0)-\theta(x)}$, where
$\theta(x):=\frac{i}{2}{\varphi(x,0)}\sigma_3$. Consequently, we
obtain
\begin{equation}
({\cal U}_+)^{-1}{{\cal U}_+}'(x) =-{\frac{i}{2}}\lambda\cdot
V_0\cdot {\begin{pmatrix}
0 & 1\\
1 & 0
\end{pmatrix}}\cdot V_0^{-1}=-\lambda\cdot
{\frac{i}{2}}\cdot\begin{pmatrix}
0 & e^{i(\varphi(0,0)-\varphi(x,0))}\\
e^{-i(\varphi(0,0)-\varphi(x,0))} & 0
\end{pmatrix}
\end{equation}
\begin{Def}
We define the {\it normalized potentials} $\eta^x$ and $\eta^y$
via the following
\begin{equation}
({\cal U}_+)^{-1}\cdot{\cal U_+}'(x) \d x :=-\lambda\cdot\eta^x,
\end{equation}
\begin{equation}
({\cal U}_-)^{-1}\cdot{\cal U_-}'(y) \d y
:=-\lambda^{-1}\cdot\eta^y,
\end{equation}
\end{Def}
Clearly, they represent $su(2)$-valued forms in $x$, respectively
$y$. Using the theorem we just proved, we obtain the form of the
normalized $x$-potential $\eta^x$:

\begin{equation}
\eta^x={\frac{i}{2}}\begin{pmatrix}
0 & e^{i(\varphi(0,0)-\varphi(x,0))}\\
e^{-i(\varphi(0,0)-\varphi(x,0))} & 0
\end{pmatrix} \d x
\end{equation}
By a completely analogous reasoning (the second part of the proof
we left to the reader), we obtain the matrix $W_0=I$ and then the
expression of the normalized $y$-potential:
\begin{equation}
\eta^y=-\frac{i}{2}\begin{pmatrix}
0 & e^{-i\varphi(0,y)}\\
e^{i\varphi(0,y)} & 0
\end{pmatrix} \d y
\end{equation}
 Note that the normalized potentials $\eta^x$ and $\eta^y$ are
completely determined by the restrictions of $\varphi$ to the axes
of coordinates. Since $\varphi (x,y)$ is invariant under
Lie-Lorentz transformations, these potentials correspond uniquely
to each (weakly regular) associate family of surfaces with Gauss
curvature $-1$.

\

Considering normalized potentials is actually equivalent to giving
a Goursat problem for the sine-Gordon hyperbolic system. In the
next paragraph, we will use the loop group splitting techniques in
order to solve this initial value problem, starting from given
normalized potentials.

\

{\bf Gauging the frame and its effect on potentials}
\begin{Def} Consider a normalized frame $\cal U$. For a rotation of
smooth angle function $\theta(x,y)$ around $e_3$, we call {\it
gauged frame} the matrix
$$ \hat{\cal U}={\cal R}_0^{-1}\cdot{\cal U}\cdot
{\cal R},$$ where ${\cal R}_0:={\cal R}(0,0)$.
\end{Def}

\begin{Def}
We define the potentials of the gauged frame $\hat{\cal U}$,
$\hat{\eta^x}$ and $\hat{\eta^y}$, by
\begin{equation}
(\hat{\cal U}_+)^{-1}\cdot\hat{\cal U_+}'(x) \d x
:=-\lambda\cdot{\hat{\eta^x}},
\end{equation}
\begin{equation}
(\hat{\cal U}_-)^{-1}\cdot\hat{\cal U_-}'(y) \d y
:=-\lambda^{-1}\cdot{\hat{\eta^y}},
\end{equation}
where
\begin{equation}
\hat{\cal U}={\hat{\cal U}}_+{\hat V_-}=\hat{\cal U}_-{\hat V_+}
\end{equation}
represent the Birkhoff splittings of the gauged frame $\hat{\cal
U}$.
\end{Def}
\begin{Prop} For a normalized frame ${\cal U}$  and its
gauge-transformed $\hat{\cal U}$, the corresponding potentials
satisfy the relations
\begin{equation}
\hat{\eta}^x={\cal R}^{-1}_0\cdot{\eta}^x\cdot {\cal R}_0,\qquad
\hat{\eta}^y={\cal R}^{-1}_0\cdot{\eta}^y\cdot {\cal R}_0.
\end{equation}
\end{Prop}

\noindent{\it Proof.} A completely straight-forward computation,
based on easy matrix manipulations and the uniqueness of the
splittings yield our formulas.

Now recall the explicit formulas (23) and (24) of the normalized
potentials $\eta^x$ and $\eta^y$, respectively. The asymmetry in
the expressions came from ``normalizing" the original orthonormal
potential $F$, that is, rotating it by the angle
$\frac{\varphi(x,y)}{2}$. In order to correct that, we have to
gauge the frame appropriately, that is rotate it ``back" with the
angle $-\frac{\varphi(x,y)}{2}$, while making sure that the
initial condition ${\cal U}(0,0,\lambda)=I$ is still satisfied.

\begin{Prop}
By gauging the normalized extended frame $\cal U$ via the rotation
$\cal R$ of angle $\theta:=-\varphi(x,y)/2$, we obtain, modulo a
constant rotation, the original orthonormal frame $\hat{\cal U}=
F=(e_1, e_2, N)=F(x,y,1)$ and its extension $F(x,y,\lambda)$ via
coordinate transformation. The potentials that correspond to the
frame $F$ are
\begin{equation}
\tilde\eta^x={\cal R}^{-1}_0\cdot{\eta}^x\cdot {\cal R}_0,\qquad
\tilde\eta^y={\cal R}^{-1}_0\cdot{\eta}^y\cdot {\cal R}_0.
\end{equation}
\end{Prop}

{\it Proof}. Based on the previous proposition, the proof is
straight-forward. \ Let us consider the normalized frame ${\cal
U}$, whose gauge correspondent is $\hat{\cal U}=F$. The potentials
are linked via the formula above, where ${\cal R}_0$ represent the
specific rotation of constant angle
${\theta(0,0)}=-\frac{\varphi(0,0)}{2}$.

Consequently, we obtain the potentials corresponding to the
orthonormal frame $F$. Denoting $\varphi_0:=\varphi(0,0)$, the
potentials corresponding to the frame $F$ are given by

\begin{equation}
\tilde\eta^x=\frac{i}{2}\begin{pmatrix}
0&e^{-i(\varphi(x,0)-{\varphi_0})}\\
e^{i(\varphi(x,0)-{\varphi_0})}&0\end{pmatrix}\d x;\qquad
\tilde\eta^y=-\frac{i}{2}\begin{pmatrix}
0&e^{-i(\varphi(0,y)-{\varphi_0})}\\
e^{i(\varphi(0,y)-{\varphi_0})}&0\end{pmatrix}\d y.
\end{equation}

Remark the symmetry of the two potentials of the frame $F$. This
is an advantage over the potentials corresponding to the
normalized frame $\cal U$.

These symmetric, ``de-normalized'', potentials are of a simpler,
more general form that we can use for the unconstrained pair of
type Weierstrass.

Note that at the origin $x=y=0$, the two potentials equal
$i{\sigma_1}/2$ and $-i{\sigma_1}/2$, respectively.

\

{\bf Constructing pseudospherical surfaces from given potentials}

We now introduce symmetric potentials $\xi^x$ and $\xi^y$ of a
general form. We will show that there is a 1-1 correspondence
between these potentials and associated families of
pseudospherical immersions.

\begin{Def}
Let $\alpha:D^x=\{x|(x,0)\in D\}\rightarrow\Bbb R$,
$\beta:D^y=\{y|(0,y)\in D\}\rightarrow\Bbb R$ be smooth functions,
such that $\alpha(0)=\beta(0)$. Let
\begin{equation}
\xi^x=\frac{i}{2}\begin{pmatrix}
0&e^{-i(\alpha(x)-{\alpha(0)})}\\
e^{i(\alpha(x)-{\alpha(0)})}&0\end{pmatrix}\d x;\qquad
\xi^y=-\frac{i}{2}\begin{pmatrix}
0&e^{-i(\beta(y)-{\beta(0)})}\\
e^{i(\beta(y)-{\beta(0)})}&0\end{pmatrix}\d y.
\end{equation}
We call $\xi^x$ and $\xi^y$ {\it symmetric potentials}. We will
use the same notations and terminology for their 3 x 3
correspondents.
\end{Def}

We are now ready to prove the following:

\begin{Th}
Let ${\hat{\cal U}_+}(y,\lambda)\in{\tilde\Lambda}^*_-{\rm SU}(2)$
and ${\hat{\cal U}_-}(x,\lambda)\in{\tilde\Lambda}^*_+{\rm SU}(2)$
be the respective solutions of the following initial value
problems:

\begin{equation}
\begin{cases}
({\cal\hat{U}_+})^{-1}{\cal{\hat{U}_+}}'(x)\d x=-\lambda{\xi}^x,\\
{\cal\hat{U}_+}(x=0)=I,
\end{cases}
\end{equation}

\begin{equation}
\begin{cases}
({\cal\hat{U}_-})^{-1}{\cal{\hat{U}_-}}'(y)\d y=-\lambda^{-1}{\xi}^y,\\
{\cal\hat{U}_-}(y=0)=I,
\end{cases}
\end{equation}
where $\xi^x$ and $\xi^y$ are given by (31). Consider the set
$${\tilde D}:=\{(x,y)\in {D^x\times D^y}\text{ ; }{\hat{\cal
U}_-}(y)\cdot{\hat{\cal U}_+}(x)\in {\tilde\Lambda}^*_-{\rm
SU}(2)\cdot{\tilde\Lambda}^*_+{\rm SU}(2)\}.$$

In $\tilde D$, we perform the Birkhoff splitting

\begin{equation}
{\hat{\cal U}_-}^{-1}(y)\cdot{\hat{\cal U}_+}(x)={\hat
V_+}(x,y)\cdot {\hat V_-}^{-1}(x,y),
\end{equation}

where $\hat V_+\in {\tilde\Lambda}^*_+{\rm SU}(2)$ and $\hat
V_-\in {\tilde\Lambda}^*_-{\rm SU}(2)$

Let

\begin{equation}
\hat{\cal U}:={\hat{\cal U}}_-\hat V_+ ={\hat{\cal U}}_+ \hat V_-
\end{equation}

Then, $\hat{\cal U}$ represents the `orthonormal frame' $F$ of an
associated family of pseudospherical surfaces in Tchebychev net,
whose Tchebychev angle $\varphi(x,y)$ verifies the conditions
$\varphi(x,0)=\alpha(x)$ and $\varphi(0,y)=\beta(y)$.
\end{Th}

\noindent{\it Proof.} Proposition 1 shows the existence and
uniqueness of a solution $\varphi$ to the initial value problem
$\varphi_{xy}=\sin
{\varphi},\varphi(x,0)=\alpha(x),\varphi(0,y)=\beta(y)$. Let
$\hat{\cal U}=F$ be the orthonormal frame corresponding to the
Tchebychev parametrization of angle $\varphi$. Formulas (30) give
the symmetric potentials $\tilde{\eta^x}$ and $\tilde{\eta^y}$
corresponding to this frame $F$, as being identical with the
symmetric potentials $\xi^x$ and $\xi^y$ assigned by (31).

In order to obtain $\varphi$ explicitely as a solution, we first
integrate (uniquely) (25) and (26), and obtain $\hat{\cal U}_+$
and $\hat{\cal U}_+$. Since $\varphi(0,0)=\alpha(0)=\beta(0)$ is
provided, so is ${\cal R}_0$. We use $\hat{\cal U}_-={\cal
R}_0^{-1}{\cal U_-}{\cal R}_0$ and $\hat{\cal U}_+={\cal
R}_0^{-1}{\cal U_+}{\cal R}_0$ to obtain ${\cal U}_+$ and ${\cal
U}_-$. Next, the Birkhoff splitting
\begin{equation}
{{\cal U}_-}^{-1}(y)\cdot{{\cal
U}_+}(x)={V_+}(x,y)\cdot{V_-}^{-1}(x,y),
\end{equation}
provides $V_+, V_-$ uniquely. Hence, the normalized frame ${\cal
U}={\cal U}_-\cdot V_+$ via formula (27), is obtained in a unique
way. We apply the Sym formula, and obtain the associated family of
immersions
\begin{equation}
\psi^{\lambda}={\frac{d}{dt}}{\cal U}^{\lambda} ({\cal U}
^{\lambda})^{-1},
\end{equation}
where $\lambda = e^t$. Finally, the map $\varphi(x,y)$ represents
the angle of this parametrization, and can be written explicitely.

\begin{Rem}
The K-Lab contains a numerical implementation of this algorithm.
Starting from two arbitrary potentials of the form (31) (i.e.,
pair of initial functions $\alpha(x)$ and $\alpha(y)$), it
computes and models the corresponding family of associated
surfaces.

\

Note that factorizations are possible only in the ``big cell",
which is an open and dense subset of the domain. The K-lab
algorithm contains an in-built numerical method that `jumps' the
singularities once they are detected, and thus allows construction
and visualization of all regular patches.
\end{Rem}

\begin{Cor}
The correspondence between the pair of symmetric potentials, and
the family of associated pseudospherical surfaces of angle
$\varphi$ is a bijection.
\end{Cor}

\noindent{\it Proof.}

Let $\Sigma$ be the map from the set of associated families of
pseudospherical surfaces in Tchebychev net into the set of all
pairs of potentials of general form (31). In essence, $\Sigma $
maps the angle $\varphi$ to the pair of potentials from (30),
which in particular are of the form (31).

On the other hand, we have a reverse procedure. Theorem 6
constructs a map from any pair of potentials (31) to a certain
family of immersions of angle $\varphi$, via the frame $\hat{\cal
U}$. We will denote this map by $\Omega$. The proof of Theorem 6
shows that the map $\Omega$ is well defined.

The construction in Theorem 6 shows that $\Sigma\circ\Omega=
{id}$, which is the same with showing that every pair of
potentials (31) is of the form (30), for a uniquely determined
angle $\varphi$ that defines a family of pseudospherical
immersions $\psi^{\lambda}$.

The uniqueness of the construction method from Theorem 6 also
shows that $\Omega \circ \Sigma = id$.

This completes the proof of the Corollary. \hbox{ }\hfill$\square$

\

{\bf Example} {\it Amsler's Surface}

In Tchebychev net parametrization, this surface corresponds to an
angle $\varphi (x,y)$ that is constant on both $x$- and $y$-axes.
For some well-known surfaces, like the pseudosphere, the
Tchebychev angle $\varphi (x,y)$ is easily written as a
trigonometric function of $x$ and $y$. This is not the case for
the Amsler surface. On the other hand, we can rewrite the
sine-Gordon equation in a very simple form ([Me, St, 2]): Let $t
:= x y$ with $(x,y)\in D=\Bbb R^2$. If we express $\varphi (x,y) =
h(xy)$, with $h:{\Bbb R} \to (0,\pi)$ a differentiable function,
then for Amsler surfaces, the sine-Gordon equation is written as
the second order differential equation
$$
t h''(t) + h'(t) = \sin (h(t)).$$ A change of function $w=e^{i
\psi}$ transforms the above equation into the so-called third
Painleve equation. Since $\varphi (x,y)$ is smooth, a
straight-forward calculation yields
$$ \varphi (0,0)=\varphi (x,0)=\varphi (0,y):=\varphi_{0}$$ for
every pair $(x,y)\in D$.
Amsler ([Ams]) showed that the solution $\varphi (x,y) = h(xy)$
oscillates near $\pi$ when $t>0$ and near $0$ when $t<0$. He also
proved that the surface has two cuspidal edges corresponding to
$\varphi=0$ and $\varphi=\pi$, respectively.

We note the two straight-lines contained in the Amsler surface,
corresponding to $x=0$ and $y=0$. As an obvious consequence of the
angle being constant along the axes, the symmetric potentials (50)
of the Amsler surface can be written as
\begin{align}
\tilde{\eta^x}&=\frac{i}{2}\begin{pmatrix}
0&1\\
1&0\end{pmatrix}\d x;
\quad\tilde{\eta^y}=-\frac{i}{2}\begin{pmatrix}
0&1\\
1&0\end{pmatrix}\d y.
\end{align}
For an interactive visualization of Amsler surfaces obtained using
the generalized Weierstrass representation (60, 61) and
computational loop-group splittings, see

http://www.gang.umass.edu/gallery/k/kgallery0201.html.

\

{\bf Acknowledgement} The author would like to express her
gratitude to her PhD advisor, Josef Dorfmeister, for his constant
support and critical comments on the present work, which is
primarily based on her dissertation. Special recognition and
thanks also go to Alexander Bobenko, Franz Pedit, and Ivan
Sterling for their outstanding and inspiring works in this area,
as well as for their suggestions and remarks. A word of
appreciation goes to Nick Schmitt for creating the K-Lab
(available at http://www.gang.umass.edu/software); his
computational work made the loop-group based construction and
visualization possible, for weakly regular pseudospherical
surfaces, as well as other important classes of surfaces. Special
thanks go to the referee of the first version of this report, for
his expert advice, excellent critical comments and many
suggestions.

\bibliographystyle{alpha}

\end{document}